\documentclass{amsart}
\usepackage{header}

\begin{document}
\title[Effective Cones of Moduli Spaces]{Effective
cones of quotients of moduli spaces \\ of stable $n$-pointed curves of
genus zero.}

\author[W. F. Rulla]{William F. Rulla}
\address{Department of Mathematics \\
	University of Georgia \\
	Athens, GA 30602}
\email{rulla@math.uga.edu}
\date{November 19, 2003}

\thanks{This paper is a product of a VIGRE seminar on $\M_{0,n}$
conducted by V. Alexeev at the University of Georgia, Athens, during
the Spring of 2002.  Thanks to S. Keel for posing the question
motivating the paper, and to him, R. Varley, and 
E. Izadi for help and advice.  Thanks also to the referee for many valuable comments.  
PORTA \cite{porta} was used in calculating several
examples.  Xfig was used for the figures.}

\keywords{Moduli space, rational curve, birational geometry,
classification of morphisms/rational maps}

\subjclass[2000]{Primary: 14E05, 14H10; Secondary: 14E30}

\begin{abstract}
Let $X_n := \M_{0,n}$, the moduli space of $n$-pointed stable genus
zero curves, and let $X_{n,m}$ be the quotient of $X_n$ by the action
of $\s_{n-m}$ on the last $n-m$ marked points.  The cones of effective
divisors $\n(X_{n,m})$, $m = 0,1,2,$ are calculated.  Using this, upper
bounds for the cones $Mov(X_{n,m})$ generated by divisors with moving
linear systems are calculated, $m = 0,1$, along with the induced
bounds on the cones of ample divisors of $\M_g$ and $\M_{g,1}$.  As an
application, the cone $\n(\M_{2,1})$ is analyzed in detail.
\end{abstract}
\maketitle

\section{Introduction}

In studying an object $X$ in the category of projective varieties, it
is useful to know to which other objects $X$ admits surjective
morphisms, or more generally, dominant rational maps.  Thus, of
interest are the sets of line bundles which have basepoint free, or
more generally, non-empty, linear series.  The closures $Nef^1(X)$ and
$\n(X)$ of the corresponding cones in the real vector space $N^1(X)$
can sometimes be determined with nothing more than basic intersection
theory.

The intermediate cone $Nef^1(X) \subseteq Mov(X)\subseteq \n(X)$,
generated by all linear series $|M|$ which have base loci $\bl |M|$ of
codimension $\ge 2$, is of greater interest than $\n(X)$ in studying
rational maps $\varphi_{|D|}:X \dasharrow Z$, since it is natural to
extend a map $\varphi_{|D|}$ over a Cartier divisor $E \in \bl |D|$ by
considering $|M := D-E|$ instead.  $Mov(X)$ similarly has an upper
bound $Nem(X)$ which is potentially calculable via basic intersection
theory.  
The purpose of this paper is to introduce the ``nem'' cone and some of
its properties, and calculate it for certain finite quotients of the
moduli spaces $\M_{0,n}$ of stable $n$-pointed curves of genus zero.

Until recently it was conjectured that the cones $\n(\M_{0,n})$ were
generated by the irreducible components of the loci of singular curves
$\Delta \subseteq \M_{0,n}$.  This is now known to be false for $n\ge
6$ (\cite{HasTsch02} or \cite{Vermeire01}), but the analogous statement
is true for the quotients $X_{n,0} := \M_{0,n}/\s_n$ by the natural
action of the symmetric group. 
A natural question (S. Keel) is: for which $m$ does the analogous
statement hold for the quotients $X_{n,m} := \M_{0,n}/\s_{n-m}$, where
$\s_{n-m}$ acts on the last $n-m$ marked points?  This paper provides
the answer:  for $m\le 2$ or $n \le 5$ (Propositions \ref{P:counter},
\ref{P:effconeXn1}, and \ref{P:effconeXn2}).  Also included is a
systematic treatment and extension of some of the ad hoc calculations
of \cite{Ru01}, relating moduli spaces $\M_{g,k}$ to various
$X_{n,m}$.  It is the author's intent to use these results for ongoing
work on $\M_3$ and other spaces.

The paper is structured as follows:  \S\ref{S:notation} establishes
notation, definitions, and basic properties of $Nem(X)$.

In \S\ref{S:NSspaces} bases for the real vector spaces $N^1(X_{n,m})$
($0 \le m \le 3$) are calculated.

\S\ref{S:counterexamples} pulls back the counterexample of
\cite{HasTsch02} and \cite{Vermeire01} in $X_{6,3}$ to all $X_{n,m}$,
$n \ge 6$ and $m \ge 3$.

\S\S\ref{S:effconeXn1}-\ref{S:NemXn1} consist of
intersection-theoretical calculations of effective and nem cones and
examples. 

In \S\ref{S:relations} the calculations of $Nem(X_{2g+2,0})$ and
$Nem(X_{2g+3,1})$ are used to impose bounds on the cones of ample divisors of
$\M_g$ and $\M_{g,1}$, respectively, via surjections of these spaces
onto the corresponding loci of hyperelliptic curves.

As an application, in \S\ref{S:M21} the cone of effective divisors
$\n(\M_{2,1})$ of the moduli space $\M_{2,1}$ of stable genus 2
curves with marked point is studied in detail.

\section{Notation, definitions, and functorial properties}\label{S:notation}

\begin{notation}
Let $X$ be a projective variety, over a(n algebraically closed) field $k$.  The characteristic $c = \CHAR k$ will be arbitrary, except for the restrictions $c \ne 2$ in \S\S\ref{S:relations} and \ref{S:M21}, and $c=0$ for the applications of Mori theory in \S\ref{S:M21}. 
\begin{enumerate}

\item The description "$\Q$-factorial" will imply normality. 

\item The symbol ``$\equiv$'' will be used to denote numerical equivalence
(of Cartier divisors or one-cycles).  

\item $N_1(X)$ is the group of one-cycles modulo numerical
equivalence, with real coefficients \cite[p. 122]{kollar}. 

\item $N^1(X)$ is the group of Cartier divisors modulo numerical
equivalence, with real coefficients \cite[p. 123]{kollar}.

\item The cone of effective divisors $\n(X)$ is the closure of the
cone in $N^1(X)$ generated by (classes of) line bundles having non-zero
global sections.

\item A $\Q$-Cartier divisor $M\in \n(X)$ will be called ``moving'' if
there is a positive integer $n$ such that ($nM$ is Cartier and) the
base locus of the linear series $|nM|$ is of codimension $\ge 2$.  The
cone of moving divisors $Mov(X)$ is the closure of the subcone of
$\n(X)$ generated by moving divisors.

\item The cone of nef divisors $Nef^1(X)$ is the closure of the
subcone of $\n(X)$ generated by line bundles having basepoint free
linear systems.  Following Miles Reid, we think of the word ``nef'' as
an acronym for ``numerically eventually free.''

\item The cone of nef curves $Nef_1(X)$ is the cone in $N_1(X)$ dual
to $\n(X)$.
\end{enumerate}
\end{notation}

Following Miles Reid, 
\begin{definition}
Let $X$ be a projective variety.  A $\Q$-Cartier divisor $M\in \n(X)$
will be called ``nem'' (for ``numerically eventually moving'') if for
every prime divisor $D \subseteq X$, $M|_D\in \n(D)$.  The cone of nem
divisors $Nem(X)$ is the closure of the cone in $\n(X)$ generated by
nem $\Q$-Cartier divisors.
\end{definition}

\begin{remarks}
\begin{enumerate}
\item Clearly the nem cone is an upper bound on the moving cone:
\[
	Mov(X) \subseteq Nem(X).
\]

\item For surfaces, the nef, moving, and nem cones are the same.



\item If $D \equiv \sum \alpha_i D_i$ where $\alpha_i > 0$ and the
$D_i$ are prime divisors, then $D \in Nem(X)$ iff $D|_{D_i} \in
\n(D_i)$ for each $i$.

\item If $M \in Nem(X)$ is a prime divisor and $C \in Nef_1(M)$, then
$C \in Nef_1(X)$.
\end{enumerate}
\end{remarks}

We prove some functorial properties. 
\begin{proposition}
If $f: X\to Y$ is a dominant, finite morphism of projective varieties,
then $f^* Nem(Y) \subseteq Nem(X)$.  Thus also 
\[
Mov(Y)\subseteq Nem(Y)\subseteq f_*Nem(X).
\]
\end{proposition}
\begin{proof}
Let $N \in Nem(Y)$ be an effective Cartier divisor.  Pick an
irreducible component $M_0$ of $f^*N$ and let $N_0 := f(M_0)$. 
Then $(f^*N)|_{M_0} = (f|_{M_0})^*(N|_{N_0}) \in \n(M_0)$.   
\end{proof}

%

\begin{lemma}
If $f:X \to Y$ is a dominant, generically finite morphism of 
$\Q$-factorial, 
projective varieties, 
then $Mov(Y) \subseteq f_*Mov(X)$.  
\end{lemma}
\begin{proof}
If $N \in Mov(Y)$ is a Cartier divisor with $|N|$ having base locus of
codimension $\ge 2$, then every divisorial component of the base locus
of $|f^*N|$ must be exceptional for $f$.  We may then produce a moving
line bundle $\til N$, mapping to (a multiple of) $N$, by subtracting
these components.
%
\end{proof}

\begin{lemma}\label{L:f_*Mmoves}
Let $f: X \to Y$ be a birational morphism of projective varieties with
$Y$ normal.  If a Cartier $M$ moves on $X$ and $f_*M$ is also Cartier,
then $f_*M$ moves on $Y$.
\end{lemma}

\begin{proof}
Because $f$ has connected fibres, any divisorial component $B$ in the
base locus of $|f_*M|$ would have strict transform contained in the
base locus of $|M|$.
%
%
\end{proof}

%

\begin{proposition}\label{P:functorial}
Let $f:X \to Y$ be a birational morphism of projective, $\Q$-factorial
varieties.  Then 
\[
Mov(Y)= f_*Mov(X)\subseteq f_*Nem(X).
\]
\end{proposition}

\begin{remark}
Proposition~\ref{P:functorial} holds when $f$ is only a ``birational
contraction'' \cite[Def. 1.0]{HuKe}, i.e. a birational map such that
the exceptional locus of $f^{-1}$ is of codimension $\ge 2$.
\end{remark} 


%

\begin{definitions}
\begin{enumerate}
\item For any $0 \le m \le n$, define $X_{n,m} := \M_{0,n}/\s_{n-m}$,
where $\s_{n-m}$ acts on the last $n-m$ marked points.  We abbreviate 
$X_{n,n} = X_{n,n-1}$ by $X_n$.  
\item Consider the quotient map $\pi:X_n \to X_{n,m}$, which makes
indistinguishable the last $n-m$ marked points.  We call the first $m$
marked points of an element of $X_{n,m}$ the \textit{distinguished
marked points}. 
\item We call the collection of images of boundary divisors
of $X_n$ the \textit{boundary divisors of $X_{n,m}$}.  We denote them
by $D^i_T$, where $i$ refers to the total number of marked points on
one component of a general member of the locus ($2\le i \le n-2$) and $T$ is the subset
of $\{1,...,m\}$ labeling the distinguished marked points
on that component ($0 \le \#T \le \min\{m,n-2\}$).  
 For $X_n$, we drop the superscript $i =
\#T$.  Note that the notation has some redundancy; e.g. in $X_n$, $D_T =
D_{T^c}$.  
\end{enumerate}
\end{definitions}

\section{The groups $N^1(X_{n,m})$}\label{S:NSspaces}
We determine bases for the groups $N^1(X_{n,m})$ for $n \ge 5$ and $0
\le m \le 3$.  From \cite{Keel92}, the spaces $N^1(\M_{0,n})$ are
spanned by the boundary classes $D_T$, subject only to the relations
($D_{T} = D_{T^c}$ and)
\[
	\sum_{\substack{i,j \in T\\ k,l\notin T}}D_T \equiv
\sum_{\substack{i,k \in T\\ j,l \notin T}} D_T, \qquad \text{for
distinct $i,j,k,l \in \{1,...,n\}$.}
\]

\begin{lemma}
The $\binom{n}{2}-n$ relations $(3 \le i < j \le n, \:\: 4 \le k \le n)$ 
\begin{equation}\label{E:relations}
\sum_{\substack{1,2\in T\\ i,j \notin T}} D_T \equiv
\sum_{\substack{1,i\in T\\ 2,j\notin T}} D_T 
\qquad \text{and} \qquad 
	\sum_{\substack{1,3\in T\\ 2,k \notin T}}D_T \equiv
\sum_{\substack{1,k\in T\\ 2,3 \notin T}} D_T 
\end{equation}
form a basis for the kernel of the vector space homomorphism
\[
	\bigoplus \R \cdot D_T \to N^1(\M_{0,n}), 
\]
where the direct sum is over the set of all boundary divisors $D_T$ of 
$\M_{0,n}$.  
\end{lemma}

\begin{proof}
The dimension of the kernel is $\binom{n}{2}-n$ \cite{Keel92}; suppose
there is a relation
\[
	\sum_{3 \le i < j \le n} a_{ij}\left( \sum_{\substack{1,2\in T\\ 
i,j \notin T}}D_T - \sum_{\substack{1,i\in T\\ 2,j \notin T}} D_T \right) +
\sum_{k=4}^n b_k \left( \sum_{\substack{1,3\in T\\ 2,k \notin T}} D_T -
\sum_{\substack{1,k\in T\\ 2,3\notin T}} D_T \right) \equiv 0.  
\]
Each of the terms $D_{\{1,2,...,n\} \backslash\{i,j\}}$ appears
exactly once, having $a_{ij}$ as its coefficient, so $a_{ij} = 0$.
In the remaining sum, each $D_{\{1,2,...,.n\}\backslash \{2,k\}}$
appears exactly once, with coefficient $b_k$, so $b_k = 0$.
\end{proof}

\begin{corollary}\label{C:N^1spaces}
\begin{itemize} 
\item 
The boundary divisors of $X_{n,0}$ and $X_{n,1}$ are a basis for their
respective groups $N^1$.  
\item
The boundary divisors of $X_{n,2}$ are subject to the single relation
\[
	\sum_{i=2}^{n-2}(n-i)(n-i-1)D^i_{12} \equiv
\sum_{i=2}^{n-2}(i-1)(n-i-1) D_1^i.  
\]
\item
The boundary divisors of $X_{n,3}$ are subject only to the following
relations (which are independent for $n \ge 5$):  
\begin{align*}
\sum_{i=2}^{n-2}(n-i-1) D^i_{12} &\equiv \sum_{i=2}^{n-2}(n-i-1)D^i_{13}
\equiv \sum_{i=2}^{n-2}(n-i-1)D^i_{23} \qquad \text{and}\\
\sum_{i=3}^{n-2}(n-i)(n-i-1) &D^i_{123} +
\sum_{i=2}^{n-3}(n-i-1)(n-i-2)D^i_{12} \\
&\equiv \sum_{i=3}^{n-2}(i-2)(n-i-1) D^i_{13} +
	\sum_{i=2}^{n-3}(i-1)(n-i-2)D^i_1.   
\end{align*}
\end{itemize}
\end{corollary}
\begin{proof}
These statements follow from pushing down the relations
(\ref{E:relations}) to the appropriate $N^1(X_{n,k})$.  
\end{proof}

Corollary~\ref{C:N^1spaces} is summarized in Table \ref{Ta:N^1spaces}.  
\begin{table}[!ht]
\begin{center}
\begin{tabular}{|| l | c | c ||}
\hline
        & Picard number & \# of boundaries \\ \hline
$X_{n,0}$ & $\lfloor \frac{n}{2}\rfloor-1$& $\lfloor \frac
	{n}{2}\rfloor -1$ \\ \hline 
$X_{n,1}$ 	& $n-3$ 	& $n-3$ \\ \hline
$X_{n,2}$ 	& $2n-7$	& $2n-6$ \\ \hline 
$X_{n,3}$ 	& $4n-16$	& $4n-13$ \\ \hline
\end{tabular}
\caption{$N^1(X_{n,k}), \: n\ge 5$}
\label{Ta:N^1spaces}
\end{center}
\end{table}

\section{Counterexamples for $X_{n,m}$, $n \ge 6$ and $m \ge 3$}
\label{S:counterexamples} 
\begin{observation} 
For $n=6$, $\n(X_{6,3})$ is not generated by boundary classes.
\end{observation}

\begin{proof}
From Corollary~\ref{C:N^1spaces}, $\rho(X_{6,3}) = 8$, an ordered
basis is given by ($D^2$, $D^2_1$, $D^2_2$, $D^2_{3}$, $D^3$,
$D^3_{1}$, $D^3_2$, $D^3_3$), and in terms of this basis,
\begin{align}\label{E:vec1}
3D^2_{12} &\equiv   (-1,1,1,0,-3,1,1,-1),\\\notag 
3D^2_{13} &\equiv   (-1,1,0,1,-3,1,-1,1), \\\notag
\text{and} \qquad 3D^2_{23} &\equiv   (-1,0,1,1,-3,-1,1,1).\notag
\end{align}
From \cite{HasTsch02} or \cite{Vermeire01} the class 
\begin{multline*}
F_{\tau} \equiv D_{36} + D_{46} + D_{56} - D_{16} - D_{26}\\ 
- D_{136} - D_{146} - D_{236} - D_{246} + D_{346} + D_{356} +
D_{12} \qquad 
\end{multline*}
lies in the effective cone of $X_6$, and (using (\ref{E:vec1}))
pushes down under the quotient morphism $X_6 \to X_{6,3}$ (making
indistinguishable the last three marked points) to
\begin{align*}
	\pi_*F_{\tau} \equiv&
2D^2_3 + D^2 + D^2 - 2D^2_1 - 2D^2_2 \\
& \qquad - 2D^3_2 - 2D^3_1 - 2D^3_1 
- 2D^3_2 + 2D^3_3 + 2D^3_3 + 6D^1_{12} \\ 
=& 2 (0,0,0,1,-3,-1,-1,1). 
\end{align*}
This is clearly not in the convex hull of the boundary divisors of
$X_{6,3}$.  
\end{proof}

For $n \ge 6$, a basis for $N^1(X_{n,3})$ is given by $\{\text{boundary
divisors}\}$ $\backslash \{D^2_{12}$, $D^2_{13}$, $D^2_{23}\}$. 
In terms of an analogous ordered basis $(D^2,D^2_1,D^2_2,D^2_3,D^3,
...)$, Corol\-lary~\ref{C:N^1spaces} gives 
\begin{align*}
(n-3)(n-4) D^2_{12} &\equiv   (-2,(n-4),(n-4),0,...),\\
(n-3)(n-4) D^2_{13} &\equiv   (-2,(n-4),0,(n-4),...),\\
(n-3)(n-4) D^2_{23} &\equiv   (-2,0,(n-4),(n-4),...),
\end{align*}
so any effective combination of these has a positive coefficient in at
least one of the second or third entries.  If $\pi: X_{6+k,3+k} \to
X_{6,3}$ forgets all the distinguished points other than the first
three, and $p:X_{6+k,3+k} \to X_{6+k,3}$ makes indistinguishable all
but the first three distinguished points, then
\[
	p_*\pi^* F_{\tau} = (0,0,0,...). 
\]
Thus the effective divisor $p_*\pi^* F_{\tau}$ is similarly not in the
cone generated by boundary divisors of $X_{n,3}$ for all $n \ge 6$.

These counterexamples lift to $X_{n,m}$ for any $m \ge 3$ via the
(finite) morphisms $X_{n,m} \to X_{n,3}$, so 

\begin{proposition}\label{P:counter}
The cone generated by boundary divisors in $N^1(X_{n,m})$ is a proper
subcone of $\n(X_{n,m})$ for all $n \ge 6$ and all $m \ge
3$. \qed
\end{proposition}

\section{The cone of effective divisors of $X_{n,1}$} \label{S:effconeXn1} 
We will use \cite[Theorem 2, p. 561]{Keel92}:
\begin{theorem}\label{T:sean}(Keel)
For any scheme $Y$ the canonical map of Chow rings 
\[
	A^*(X_n) \otimes A^*(Y) \to A^*(X_n\times Y)
\]
is an isomorphism.
\end{theorem}

\begin{corollary}\label{C:productcones}
If $V$, $W$ are finite images of some $X_n$, $X_m$,
respectively, then
\begin{align*}
	N^1(V\times W) &\cong p_1^* N^1(V)\times p_2^* N^1(W) \\
\text{and so} \qquad 	\n(V\times W) &\cong p_1^* \n(V)\times p_2^* \n(W),
\end{align*}
where $p_i$ is the $ith$ projection.  \qed
\end{corollary}

For notational simplicity, we make the following definitions for
divisors on $X_{n,m}$ ($0 \le m \le 2$):
\begin{definitions}
\begin{enumerate}

\item Set $B_k := D^k_{\{1,...,m\}}$ ($2\le k \le n-2$).  Note that
when $m=0$, $B_{n-k} = B_k$.

\item It is straightforward to show that divisorial ramification of
the quotient morphism $X_n \to X_{n,m}$ occurs only above $B_{n-2}$,
and the ramification is simple (see \cite{KeMc96} for a full
discussion).  Thus, set $b_{n-2} := \frac{1}{2} B_{n-2}$, and $b_i :=
B_i$ for all other $i$.  It is standard practice in the intersection
theory on $\M_g$ to deal with ramification at this stage
(cf. $\delta_1 = \tfrac{1}{2}\Delta_1$); using the $b$'s instead of
the $B$'s accounts for the translation between intersection
calculations in $X_n$ and $X_{n,m}$.  The notation also simplifies
formulae by allowing all the boundary classes to be treated in a
symmetric way.

\item 
For $2 \le k \le n-3$ let $C_k$ be the curve class obtained by
varying the node of a general member of $B_k$ along the component not
containing any distinguished marked point.  
\begin{figure}[!ht]
\begin{center}
	\scalebox{.3}{\includegraphics{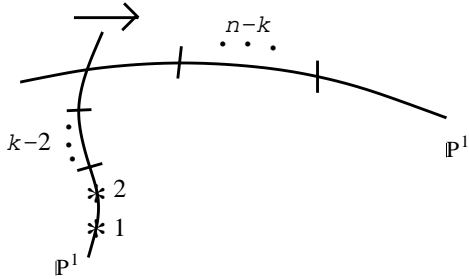}}
	\caption{The one-parameter family $C_k$. }\label{Fi:Ck}\label{Ck}
\end{center}
\end{figure}

\end{enumerate}
\end{definitions} 

\begin{notation}
Given divisor classes $d_1,...,d_k\in N^1(X)$ and a curve class $C\in
N_1(X)$,
\begin{equation}\label{E:notation}
``\:C \equiv \sum n_i\check d_i \: ''
\end{equation}
will mean $C.d_i = n_i$ for each $i$.  If the $d_i$ form part of a
basis of $N^1(X)$, $\check d_i$ will denote the dual to $d_i$ in the
dual basis, and (\ref{E:notation}) will further mean that the basis
elements whose duals do not appear intersect $C$ trivially.
\end{notation}
Now fix $m=1$.  The $\{b_i\}$ give a basis for $N^1(X_{n,1})$; we
formally define $\check b_1 = 0$.  Define $C_1$ to be the curve class
of a fibre of the morphism $X_{n,1} \to X_{n-1,0}$ forgetting the
distinguished marked point.  Then
\begin{equation}\label{E:Ck}
	C_k \equiv (n-k)\check b_{k+1} + (2-n+k)\check b_k \qquad (1
\le k \le n-3).
\end{equation}

$C_k \in Nef_1(B_k)$ for $k\ge 2$, and $C_1 \in Nef_1(X_{n,1})$.  Thus
any prime divisor $D$ having class not proportional to that of any
boundary divisor has $D.C_k \ge 0$ for all $k = 1, ..., n-3$.
Rewriting (\ref{E:Ck}) as
\[
	\check b_{k+1} \equiv \frac{1}{n-k}\bigl[ C_k +
(n-k-2)\check b_k \bigr] \qquad (1 \le k \le n-3),
\]
we see inductively that the coordinates of any such $D$ are
non-negative, so

\begin{proposition}[$n\ge 4$]\label{P:effconeXn1}
$\n(X_{n,1})$ is simplicial, generated by the boundary classes $B_i$
$(i = 2,...,n-2)$.
\end{proposition}

\section{$Nem(X_{n,0})$}\label{S:NemXn0}
We use Proposition~\ref{P:effconeXn1} to calculate $Nem(X_{n,0})$.
\begin{lemma}\label{L:effconeXn0}
$\n(X_{n,0})$ is generated by its boundary divisors.
\end{lemma}

\begin{proof}
This is shown in \cite{KeMc96}; it is also an immediate corollary of
Proposition~\ref{P:effconeXn1}. 
\end{proof}

The boundary divisors of $X_{n,0}$ are of the form $X_{l+1,1}\times
X_{n-l+1,1}$, so by Lem\-ma~\ref{L:effconeXn0} and
Theorem~\ref{T:sean} we need only consider contributions from
$X_{l+1,1}$ $(3 \le l \le n-2)$.  More generally, we will calculate the
contribution of $X_{l+1,1}$ to the inequalities determining
$Nem(X_{n,m})$, $m=0,1,2$.  For $m=0,1$, a basis is given by the set
of all $b_i$ (note that $b_{n-k}=b_k$ when $m=0$); we formally define
$\check b_1 :=0$.  For $m=2$ we will use $b_3,...,b_{n-2}$ as part of
a basis for $N^1(X_{n,2})$, and formally define $\check b_1 = \check
b_2 = 0$.

For $3 \le l \le n-2$, $l \le n-m$, and $0 \le m \le 2$ let 
\[
	q: X_{l+1,1} \to X_{n,m}
\]
be the inclusion given by attaching to the distinguished marked point
of every member of $X_{l+1,1}$ a fixed element of $X_{n-l+1,m+1}$ at
its highest numbered distinguished point.  One checks ($1 \le k \le
l-1$):
\[
	q_* C_k \equiv (l-k+1)\check b_{n-l+k} - (l-k-1)\check
	b_{n-l+k-1}, 
\]
giving the recursive relationship
\begin{equation}
q_* \check {b}_{k+1} = \check b_{n-l+k} - \left( \frac{l-k-1}{l-k+1}
\right) \check b_{n-l+k-1} + \left( \frac{l-k-1}{l-k+1} \right)
qm_{l+1*}\check b_k.
\end{equation}
The recursion telescopes, 
\begin{equation}\label{E:initQnContrib}
q_* \check b_{k+1} = \check b_{n-l+k} -
\frac{(l-k-1)(l-k)}{l(l-1)} \check b_{n-l},
\end{equation}
giving the non-negative half-planes $(1 \le k \le l-2, \:\: 3 \le l
\le n-2)$
\[
	l\cdot (l-1) \check b_{n-l+k} \ge (l-k-1)(l-k) \check
b_{n-l}.  
\]
The inequalities for $k=1$ generate, so (recalling for $m=2$ that
$\check b_2 := 0$)   
\begin{proposition}[$m = 0,1,2$]\label{P:QnContrib}
The contribution from the spaces $X_{l+1,1}$ to the inequalities
determining $Nem(X_{n,m})$ are
\[ 
	l \cdot \check b_{n-l+1} \ge (l-2) \cdot \check b_{n-l} \qquad 
(3 \le l \le n-2).
\]
\end{proposition}

For $m=0$, $b_{n-k} = b_k$, and Proposition~\ref{P:QnContrib} and
Lemma~\ref{L:effconeXn0} imply:

\begin{theorem}[$n \ge 6$]\label{T:NemZn}
$Nem(X_{n,0})$ is determined by the $2(\lfloor
\frac{n}{2}\rfloor -2)$ inequalities $(2 \le i \le \lfloor
\frac{n}{2}\rfloor -1)$:
\[
\left(\frac{n-i-2}{n-i} \right) \check{ b}_i \: \le \:
	\check{ b}_{i+1} \: \le \: 
		\left( \frac{i+1}{i-1} \right) \check{ b}_i.
\]
In particular, any moving linear system of $X_{n,0}$ is big.
\qed
\end{theorem}

Because of the simplicity of these inequalities, we can explicitly
describe the cone $Nem(X_{n,0})$.
\begin{corollary}
$Nem(X_{n,0})$ is generated by $2^{\lfloor \tfrac{n}{2} \rfloor
-2}$ rays, which are inductively described in terms of the basis
$\{\check { b}_i\}$ as follows. Let the first entry of a generating
vector be $1$.  Then given that the $i$th entry is $a$,
there are two choices for the $i+1$st, namely 
\[
	a \cdot \left( \frac{n-i-2}{n-i} \right) \qquad \text{and} 
\qquad a \cdot \left( \frac{i+1}{i-1} \right). 
\] \qed
\end{corollary}

\begin{remark}
For any $i \in \{2,...,\lfloor \tfrac{n}{2}\rfloor\}$, there is an
extremal ray $R_i \subseteq Nem(X_{n,0})$ which, when restricted to
each of the $B_j$, is big only on $B_i$.  One constructs a vector $v_i$
generating $R_i$ as follows:  let the left-hand inequalities of
Theorem~\ref{T:NemZn} be called ``Type A'' and the others ``Type B.''
Set the $ith$ entry of $v_i$ equal to $1$, and use the inequalities of
Type A (as equalities) to determine the entries to the right of the
$ith$, and those of Type B to determine those to the left.

The $R_i$ are interesting from the point of view of Mori theory: if
they were generated by moving linear systems having finitely generated
section rings (in particular if the $X_{n,0}$ were Mori Dream Spaces), 
they would correspond to birational maps $\varphi_i: X \dasharrow Z_i$
contracting all the $B_j$ except $B_i$; the targets $Z_i$ are
$\Q$-factorial of Picard number one (see \cite{HuKe} for background on
section rings and ``contracting" rational maps).
\end{remark}

\begin{examples}
We illustrate cross-sections of $\n(X_{n,0})$ for $6\le n \le 9$.  In
the following, we use as ordered bases $(b_2,..., b_{\lfloor
\frac{n}{2}\rfloor})$.  For comparison purposes, we have included the
calculations of $Nef^1(X_{n,0})$ done in \cite{KeMc96}.  Moving and
nem cones are the same for $n\le 6$, 
but this is not known for $n \ge 7$.

\begin{enumerate}
\item[$\mathbf{n=6}$:] $Nem(X_{6,0}) = Nef^1(X_{6,0})$ is generated by
$(2,1)$ and $(1,3)$.

\begin{figure}[!ht]
\begin{center}
	\scalebox{.35}{\includegraphics{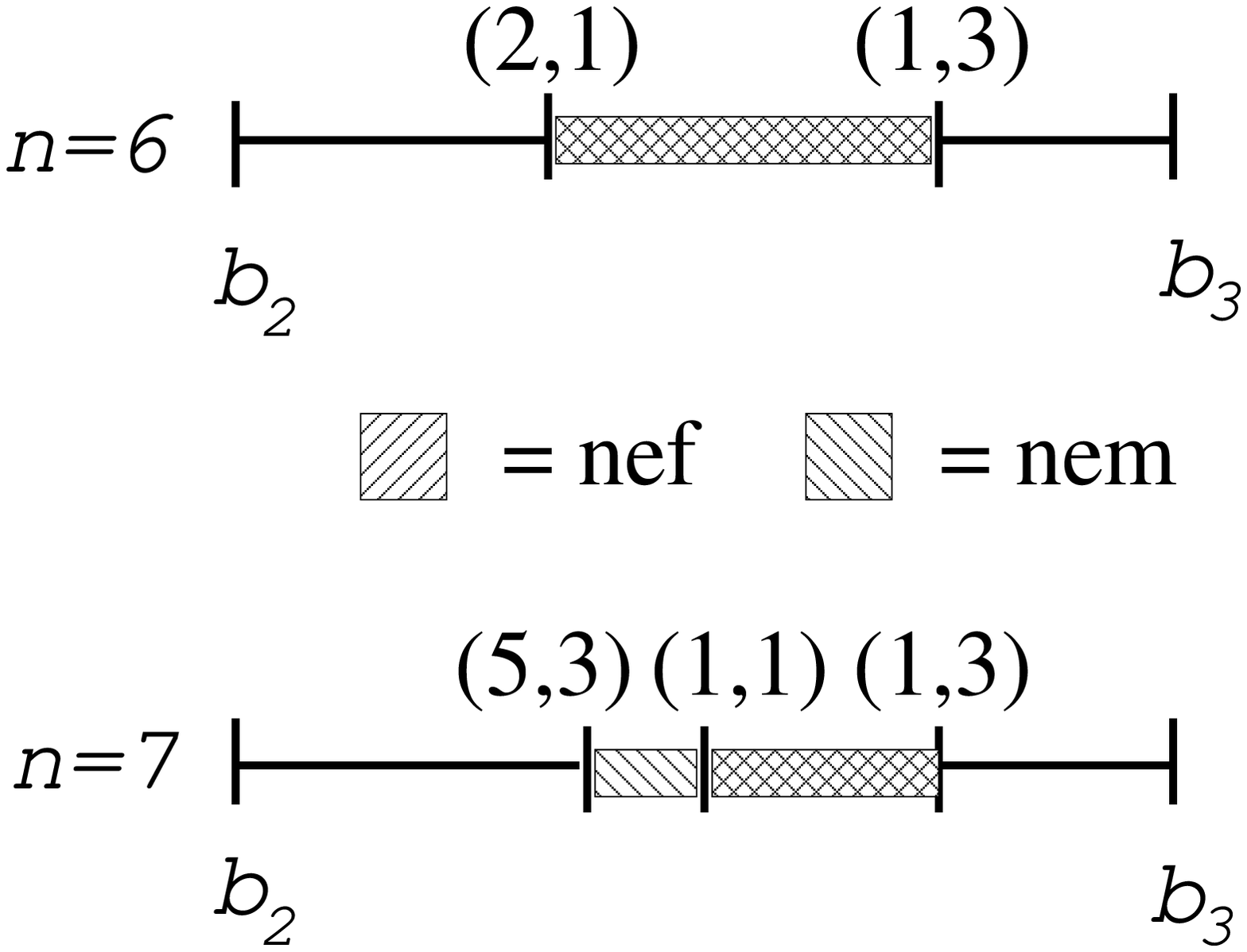}}
	\caption{Cross sections of $\n(X_{6,0})$ and $\n(X_{7,0})$. }\label{Fi:Z67}\label{Z67}
\end{center}
\end{figure}

\item[$\mathbf{n=7}$:] $Nem(X_{7,0})$ is generated by $(5,3)$ and $(1,3)$.
$Nef^1(X_{7,0})$ is generated by $(1,3)$ and $(1,1)$.

\item[$\mathbf{n=8}$:] $Nem(X_{8,0})$ is generated by 
\[
A := (3,2,4), \quad C := (1,3,6), \quad Q := (5,15,9), \quad R := (15,10,6).
\]
$Nef^1(X_{8,0})$ is spanned by $A$ and $C$ together with $P := (6,11,8)$
and $B := (2,6,5)$.

\begin{remark}
The three vertices $A$, $B$, and $C$ come from the nef cone of $\M_3$
under the natural identification of $X_8$ with the hyperelliptic locus
in $\M_3$ \cite{Ru01}.
\end{remark}

\begin{figure}[!ht]
\begin{center}
	\scalebox{.35}{\includegraphics{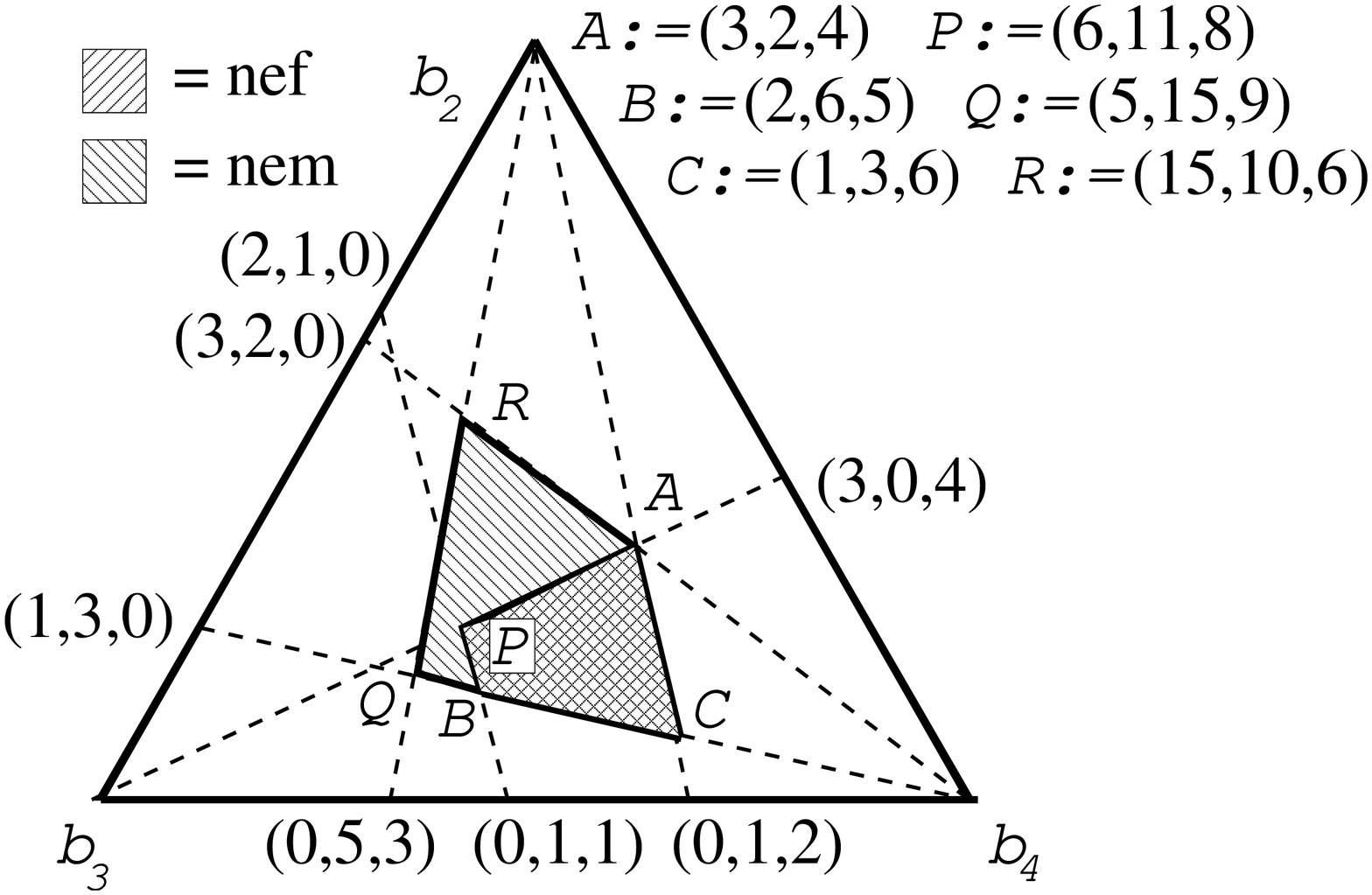}}
	\caption{A cross section of $\n(X_{8,0})$. }\label{Fi:M08tilde}\label{M08tilde}
\end{center}
\end{figure}

\item[$\mathbf{n=9}$:] $Nem(X_{9,0})$ is generated by
\[
B := (1,3,2), \quad C := (1,3,6), \quad T := (7,5,10), \quad U := (21,15,10).
\]
$Nef^1(X_{9,0})$ is spanned by $B$ and $C$ together with $A :=
(1,1,2)$ and $S := (3,3,4)$.

\begin{figure}[!ht]
\begin{center}
	\scalebox{.35}{\includegraphics{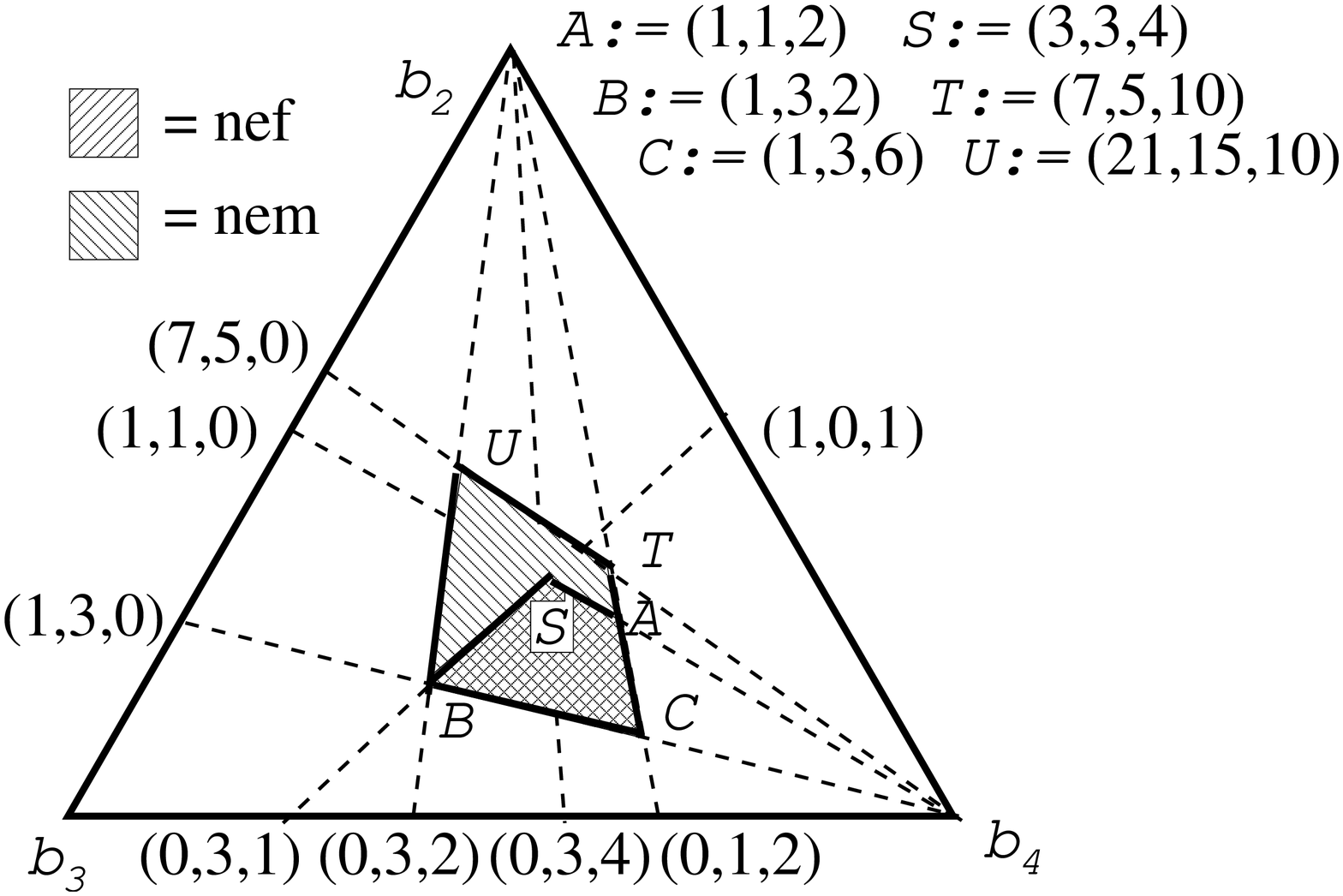}}
	\caption{A cross section of $\n(X_{9,0})$.}\label{Fi:M09tilde}\label{M09tilde}
\end{center}
\end{figure}

\end{enumerate}
\end{examples}

\section{The cone of effective divisors of $X_{n,2}$}
For notational ease, define $b^*_i := D_{\{1\}}^i$ in $X_{n,2}$ .  From
Corollary~\ref{C:N^1spaces} a basis for $N^1(X_{n,2})$ is given by the
set of all $\{b_i, b^*_i\}$ excluding $b_2$.  Recall that in
$N_1(X_{n,2})$ we formally define $\check b_1 = \check b_2 = 0$; we
also set $\check b_1^* = 0$.  

Let $D$ be any prime divisor having class not proportional to that of any
boundary. The boundaries of $X_{n,2}$ are of the form $X_{l+1,3}\times X_{n-l+1,1}$ or $X_{l+1,2}\times X_{n-l+1,2}$.  Thus the restriction to $D$ to any appropriate factor $X_{l+1,m}$ ($1\le m \le 3$) must be effective.  
Constraints coming from spaces of the form $X_{l+1,1}$ were determined
in Proposition~\ref{P:QnContrib}; in particular, $\check b_i.D \ge 0$
$(3 \le i \le n-2)$.  

Next consider constraints coming from the spaces $X_{l+1,2}$ $(2 \le l
\le n-2)$.
\begin{lemma}\label{L:EffRn} 
The subcone of $N^1(X_{l+1,2})$ generated by the classes $\{b_i,
b^*_i\}_{i=2}^{l-1}$ is determined by the inequalities $(2 \le j \le
l-1, \: \: 3 \le i \le l-1)$:
\begin{align*}
\check b^*_j &\ge 0 \qquad \text{and} \\
(j-1)(l-j)\check b_i + &(l-i+1)(l-i)\check b^*_j \ge 0.
\end{align*}
\end{lemma}
\begin{proof}
This follows from Corollary~\ref{C:N^1spaces}. 
\end{proof}


\begin{definition}
Define the curve class $C_1^*$ by varying the marked point $1$ along a
general element of $X_{l+1,2}$.  For $2 \le i \le l-2$ define $C^*_i$ by
varying the component of a general element of $b^*_i \subseteq
X_{l+1,2}$ with the marked point $1$ and $i-1$ indistinguishable marked
points along the component with the remaining marked point.
\end{definition}
Thus, $C^*_i \equiv \check b_{i+1} + (l-i) \check b^*_{i+1} + (1-l+i)
\check b^*_{i}$,  or
\[
\check b^*_{i+1} \equiv \frac{1}{l-i} \biggl( C^*_i -
\check b_{i+1} + (l-i-1) \check b^*_i \biggr) \qquad (1 \le i \le
l-2).
\]

Consider the morphism $r:X_{l+1,2} \to X_{n,2}$ which attaches a fixed
element $(C,\bar P_1, \bar P_2)$ $\in$ $X_{n-l+1, 2}$ to each element $(C,
P_1$, $P_2)$ $\in X_{l+1,2}$ by identifying $P_2$ with $\bar P_1$.
Then a calculation similar to that done in \S\ref{S:NemXn0} shows 
\begin{equation}\label{E:twototwopushforwards}
r_* \check b^*_{i+1} \equiv \check b^*_{i+1} +
\frac{i(l-i-1)}{(l-2)(l-1)} \check b_{n-l+1} - \frac{i}{l-1} \check
b^*_l \qquad (1 \le i \le l-2).
\end{equation}
Pushing forward the inequalities $\check b^*_j \ge 0$ of 
Lemma~\ref{L:EffRn} gives 
\[
	\check b^*_j + \frac{(j-1)(l-j)}{(l-2)(l-1)} \check b_{n-l+1}
	\ge \left(\frac{j-1}{l-1}\right) \check b^*_l \qquad (2 \le j
	\le l-1). 
\]
Fixing $l = n-2$ we get necessary inequalities for the prime divisor $D$:  
\begin{equation}\label{E:ineq1}
(n-4)(n-3) \check b^*_j + (j-1)(n-j-2) \check b_3 \ge (n-4)(j-1) \check
b^*_{n-2}.  
\end{equation}
If we reverse the roles of the two marked points (in the definition of
$C_i^*$ and the morphism $r_k$) we get   
\begin{equation}\label{E:ineq2}
(n-4)(n-3) \check b^*_{n-j} + (j-1)(n-j-2)  \check b_3  \ge (n-4)(j-1) 
\check b^*_2.   
\end{equation}
Replacing $j$ with $n-j$ in (\ref{E:ineq2}) gives  
\begin{equation}\label{E:ineq3} 
	(n-4)(n-3) \check b^*_j + (n-j-1)(j-2) \check b_3 \ge (n-4)(n-j-1)
\check b^*_2.
\end{equation}
Note that (\ref{E:ineq1}) - (\ref{E:ineq3}) continue to hold for
$j=n-2$, so apply for $2 \le j \le n-2$. 

Varying the point of attachment of the component with $n-3$
undistinguished marked points of a general element of $B_3 \subseteq
X_{n,m}$ along the other component gives
\begin{equation}\label{E:ineq4}
\check b^*_2 + \check b^*_{n-2} \ge \check b_3, 
\end{equation}
satisfied by any prime $D$ as above.  Then $\check b^*_j \ge 0$ $(2
\le j \le n-2)$ comes from
\[
(n-j-1) [\text{Ineq.} (\ref{E:ineq1})] 
+ (j-1) [\text{Ineq.} (\ref{E:ineq3})]
+ (n-j-1)(j-1)(n-4) [\text{Ineq.} (\ref{E:ineq4})].
\]

\begin{proposition}[$n \ge 4$]\label{P:effconeXn2}
$\n(X_{n,2})$ is generated by boundary divisors.  Further,
$Nem(X_{n,2})$ is contained in the (simplicial) cone generated by the
classes $\{b_i, b^*_j\}$ excluding $b_2$.  
\qed
\end{proposition}

\begin{example}
Consider the case $n=5$.  From Corollary~\ref{C:N^1spaces},
$\rho(X_{5,2}) = 3$ and
\[
	3b_2 \equiv b^*_2+b^*_3 - b_3.  
\]
$\n(X_{5,2})$ is generated by the four boundaries $b_2$, $b_3$,
$b_2^*$, and $b_3^*$. $Nef^1(X_{5,2})$ is the dual of this cone; since
$X_{5,2}$ is a surface, its nem, moving, and nef cones are identified
(Figure~\ref{Fi:X53}).

\begin{figure}[!ht]
\begin{center}
	\scalebox{.35}{\includegraphics{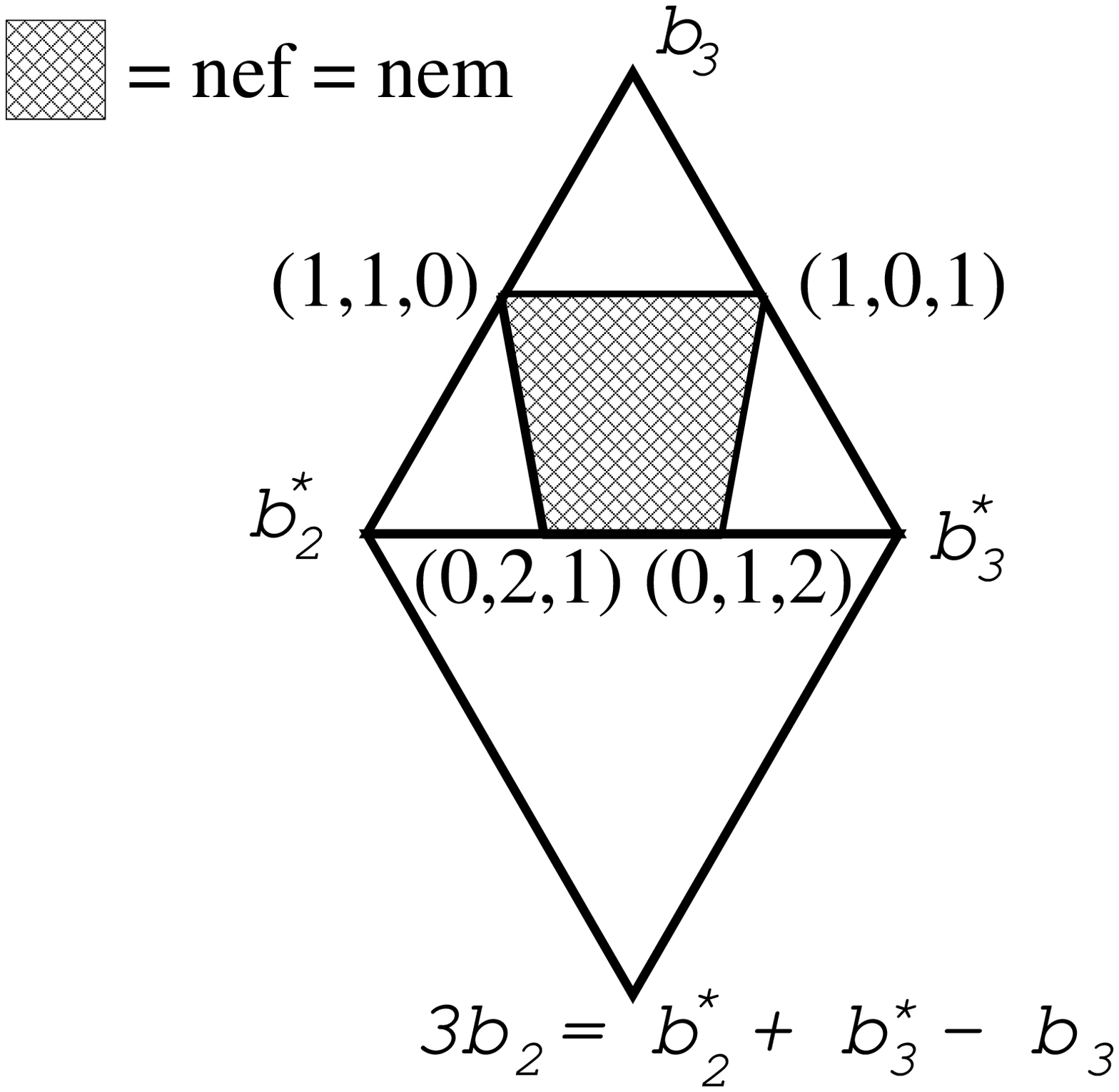}}
	\caption{A cross section of $\n(X_{5,2})$. }\label{Fi:X53}\label{X53}
\end{center}
\end{figure}
\end{example}

\section{$Nem(X_{n,1})$}\label{S:NemXn1} 

We use Proposition~\ref{P:effconeXn2} to calculate $Nem(X_{n,1})$.
The boundaries of $X_{n,1}$ are of the form $X_{l+1,2}
\times X_{n-l+1,1}$.  Proposition \ref{P:QnContrib} gives the
contributions from the $X_{l+1,1}$:  
\begin{equation}\label{E:QnContrib2}
	l \cdot \check b_{n-l+1} \ge (l-2) \cdot \check b_{n-l} \qquad
	(3 \le l \le n-2).
\end{equation}

If $s: X_{l+1,2} \to X_{n,1}$ is a corresponding attaching map,
we adapt (\ref{E:initQnContrib}) and
(\ref{E:twototwopushforwards}) to give ($2 \le i \le l-2$):
\begin{equation}\label{E:twotoonepushforwards} 
\begin{aligned}
s_*\check b_{i+1} &\equiv \check b_{n-l+i} -
\frac{(l-i-1)(l-i)}{(l-2)(l-1)} \check b_{n-l+1} \\
\text{and} \qquad s_*\check b^*_{i+1} &\equiv \check b_{i+1} +
\frac{i(l-i-1)}{(l-2)(l-1)}\check b_{n-l+1} - \frac{i}{l-1}\check
b_l.
\end{aligned}
\end{equation}
Pushing forward the inequalities determining $\n(X_{l+1,2})$
(Lemma~\ref{L:EffRn} and Proposition~\ref{P:effconeXn2}) gives 
\begin{multline}\label{E:firstTry}
(l-1)(j-1)(l-j)\check b_{n-l+i-1} + (l-1)(l-i)(l-i+1)\check b_j \\
\: \ge \: (j-1)(l-i)(l-i+1) \check b_l \qquad (2\le i,j\le l-1). 
\end{multline}

\begin{corollary}\label{C:NemXn1} 
$Nem(X_{n,1})$ is generated by the $(n-1)(n-4)/2$ inequalities
$(3\le l \le n-2, \: \: 2 \le j \le l-1)$: 
\begin{equation}\label{E:secondTry}
\begin{aligned}
l \cdot \check b_{n-l+1} &\ge (l-2) \check b_{n-l}           
\qquad \text{and} \\
(j-1)(l-j)\check b_{n-l+1} + (l-1)(l-2)\check b_j &\ge
(j-1)(l-2)\check b_l.  \hskip1in
\end{aligned}
\end{equation}
\end{corollary}
\begin{proof}
$Nem(X_{n,1})$ is by definition a subset of $\n(X_{n,1})$, so the
inequalities $\check b_i \ge 0$ must be included in its calculation.
These however are subsumed by (\ref{E:QnContrib2}) and
(\ref{E:firstTry}): let $I_{i,j,l}$ refer to the inequality 
(\ref{E:firstTry}).  If $n$ is odd, ``$\check b_2 \ge 0$'' comes from
$I_{2,2,(n+1)/2}$; if $n$ is even, ``$\check b_2 \ge 0$'' comes from
$(\tfrac{n}{2})I_{2,2,n/2} + (\tfrac{n}{2}-2)I_{2,2,n/2+1}$.  The
other inequalities ``$\check b_i \ge 0$'' then follow from
(\ref{E:QnContrib2}).

Denote the inequalities (\ref{E:QnContrib2}) by $J_l$.  For $i \ge 3$,
\[
	I_{i,j,l} = a J_{l-i+2} + bI_{i-1,j,l},  
\]
where 
\[
a = \frac{(l-1)(j-1)(l-j)}{(l-i+2)} \qquad \text{and} \qquad 
b = \frac{l-i}{l-i+2}.  
\]
Thus we need only the $J_l$ and the $I_{2,j,l}$.  
\end{proof}

We now study some properties of $Nem(X_{n,1})$.  From the inequalities
(\ref{E:QnContrib2}) note that if an element of $Nem(X_{n,1})$ has
$\check b_i = 0$, then it also has $\check b_j = 0$ for all $j <i$.
Suppose that $\check b_3 = 0$.  Then $I_{2,2,n-2}$ implies $0 \ge
\check b_{n-2}$, so that all $\check b_i = 0$.  Thus the only possible
zero coordinate of a vector $(v\ne 0) \in Nem(X_{n,1})$ is the first:
$\check b_2 (v) =0$.  This implies $v$ trivially intersects a fibre of
$\pi:  X_{n,1} \to X_{n-1,0}$, so we expect $v$ to be pulled back from
$Nem(X_{n-1,0})$.  In fact, if $\check b_2 =0$, then $I_{2,2,l}$
implies
\[
	\check b_{n-l+1} = \check b_l \qquad (3 \le l \le n-2);
\]
together with the inequalities $J_l$, this forces 
\[
	\frac{n-l-2}{n-l} \check b_l \;\le\: \check b_{l+1} \:\le\:
\frac{l}{l-2} \check b_l.  
\]
Since $\pi^* b_l = b_{l+1} + b_{n-l+1}$ (except $\pi^* b_{\frac{n}{2}}
= b_{\frac{n}{2}+1}$ for $n$ even) we get from Theorem~\ref{T:NemZn}
\begin{proposition}
	$Nem(X_{n,1}) = \pi^*Nem(X_{n-1,0}) + C$, where every element of
$C$ is big. \qed
\end{proposition}

Thus any fibration of any space $V$ isomorphic to $X_{n,1}$ in
codimension one contracts the (images of the) fibres of $\pi:X_{n,1}
\to X_{n-1,0}$.

\begin{examples}  We illustrate cross sections of $\n(X_{n,1})$ for $5
\le n \le7$.  In the following we again include the computations of
nef cones from \cite{KeMc96}.

\begin{figure}[!ht]
\begin{center}
	\scalebox{.35}{\includegraphics{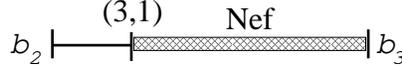}}
	\caption{A cross section of $\n(X_{5,1})$. }\label{Fi:Q5}\label{Q5}
\end{center}
\end{figure}

\begin{enumerate}

\item[$\mathbf{n =5}$:]  From Corollary~\ref{C:NemXn1}, 
$Nem(X_{5,1})$ is determined by $3 \check b_3 \ge \check b_2$ and
$\check b_2 \ge 0$.  Since $X_{5,1}$ is a surface, this also gives the
nef cone (Figure~\ref{Fi:Q5}).

\item[$\mathbf{n = 6}$:] From Corollary~\ref{C:NemXn1}, $Nem(X_{6,1})$ is
determined by 
\[
\begin{aligned}
3\check b_4 &\ge \check b_3, \\
4\check b_3 &\ge 2\check b_2,
\end{aligned}
\qquad 
\begin{aligned}
2\check b_4 + 4\check b_2 &\ge 2\check b_3, \\  
6\check b_3 + 18\check b_2  &\ge 6\check b_4,
\end{aligned}
\qquad 
\begin{aligned}
6\check b_3 + 18\check b_3  &\ge 12\check b_4.
\end{aligned}
\]
These determine the cone spanned by the vectors 
\[
(6,3,1), \quad (1,3,1), \quad (0,1,1), \quad (1,3,6), \quad (2,1,2).
\]

The nef cone is given by the span of 
\[
(6,3,1), \quad (1,3,1), \quad (0,1,1), \quad (2,1,2).
\]
\begin{remark}
$(6,3,1)$ corresponds to (i.e. is a multiple of the pull-back of) the
relative dualizing sheaf $\omega$ of $\M_{2,1}$ under the
identification of $X_{6,1}$ with the Weierstrass locus of $\M_{2,1}$
(see \S\ref{S:M21} for more information), $(2,1,2)$ corresponds to the
determinant of the Hodge bundle $\lambda$, and $(1,3,1)$ corresponds
to $12\lambda - \delta_{irr}$, where $\delta_ {irr}$ denotes the class
of the closure of the locus of irreducible elements of $\M_{2,1}$ with
a single node.
\end{remark}
\begin{figure}[!ht]
\begin{center}
	\scalebox{.35}{\includegraphics{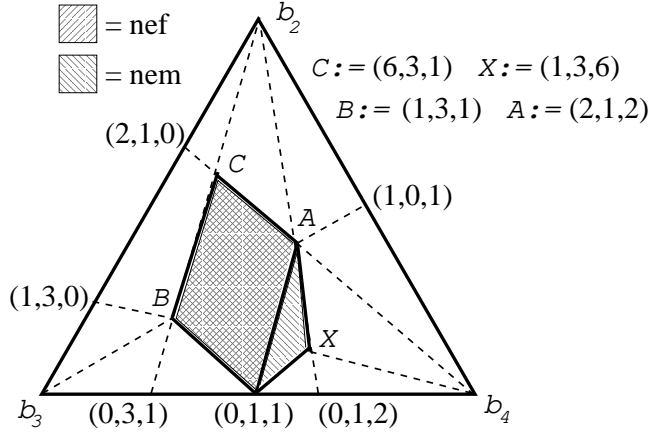}}
	\caption{A cross section of $\n(X_{6,1})$. }\label{Fi:Q6}\label{Q6}
\end{center}
\end{figure}

\item[$\mathbf{n=7}$:]  From Corollary~\ref{C:NemXn1}, $Nem(X_{7,1})$ is
determined by    
\[
\begin{aligned}
3 \check b_5 &\ge \check b_4, \\ 
4\check b_4 &\ge 2 \check b_3,
\end{aligned}
\qquad 
\begin{aligned} 
5 \check b_3 &\ge 3 \check b_2,\\
2 \check b_5 + 4 \check b_2 &\ge 2 \check b_3,
\end{aligned}\qquad  
\begin{aligned}
\check b_2 &\ge 0, \\  
\check b_3 + 4 \check b_2 &\ge \check b_5,
\end{aligned}
\qquad 
\begin{aligned}
8 \check b_3 &\ge 3 \check b_5,\\ 
\check b_3 + 4 \check b_4 &\ge 3 \check b_5.
\end{aligned}
\]
This gives $Nem(X_{7,1})$ as the convex hull of 
\begin{align*}
(5,12,36,32),& \quad (10,6,3,6),& \quad &(5,3,9,8),& \quad (5,12,6,2),& \quad (0,1,3,1),& \\
(5,12,21,32),& \quad (10,6,3,1),& \quad &(5,3,9,3),& \quad (20,12,21,32),& \quad (0,2,1,2).& 
\end{align*}
$Nef^1(X_{7,1})$ is the cone over
\[
(10,6,3,1), \quad (5,3,4,3), \quad (5,12,6,2),\quad (0,1,3,1), \quad (0,2,1,2).
\]
\begin{remarks}
\begin{enumerate}
\item
$(0,2,1,2)$ corresponds to $\lambda$ under the morphism $p:X_{7,1}\to
\M_{2,1}$ described in \S\ref{S:M21}.  Also, $(0,1,3,1)$ corresponds to
$12\lambda - \delta_{irr}$, and $(5,12,6,2)$ corresponds to $\omega$.

\item \label{I:moving}
It is easy to check that the ray generated by $(10,6,3,1)$ is the push-forward
under $\pi_7: X_7 \to X_{7,1}$ of the ray generated by 
\begin{multline*}
\qquad L_7 \equiv D_{12} +  D_{37}+ D_{47}+ D_{57}+ D_{67}+ D_{347}+
D_{357} + \\ 
D_{367}+ D_{457}+ D_{467}+ D_{567}+ D_{3457}+ D_{3467} + D_{3567}+
D_{4567},
\end{multline*}
which corresponds to a hyperplane section from $\P^4$ under the
blow-up description $f:\M_{0,7}\to \P^4$ in \cite{HasTsch02}.  $|L_7|$
contains elements missing any given fibre of $\pi_7$, so $(10,6,3,1)$
moves in $X_{7,1}$.  We use this in the proof of Proposition
\ref{P:MovM21}.
\end{enumerate}
\end{remarks}

\end{enumerate}
\end{examples} 

\section{Relations with $\M_g$ and $\M_{g,1}$}\label{S:relations} 
We use the notation of \cite{FaberAlgs} for the standard divisor classes on $\M_{g}$ and $\M_{g,1}$, except on $\M_{g,1}$ we abbreviate $\delta_{i;\{1\}}$ by $\delta_i$ ($1\le i \le g-1$), and $\omega_{\pi_1}$ by $\omega$.  It is a fact that the classes $\delta_{irr}$, $\delta_1$, ..., $\delta_{\lfloor\tfrac{g}{2}\rfloor}$ and $\lambda$ generate $N^1(\M_g)$, and the classes $\delta_{irr}$, $\delta_{1}$, ..., $\delta_{g-1}$, $\lambda$, and $\omega$ generate $N^1(\M_{g,1})$ in any characteristic \cite{DM69}, \cite{Moriwaki2}.  These classes are independent for $g \ge 3$.  

To avoid pathologies with hyperelliptic curves, we restrict to characteristic $c \ne 2$ in this and the following section.   

\subsection{The hyperelliptic locus of $\M_g$}\label{SS:Mg}
Consider the natural morphism 
\[
i: X_{2g+2,0} \to \M_g
\]
identifying the source with the hyperelliptic locus of $\M_g$, by
associating to each element of $X_{2g+2,0}$ its (unique) degree-two
admissible cover \cite[pp. 180 ff.]{Mgbook}.  For $n \ge 3$ we use the ordered basis $(\lambda, \delta_{irr}, \delta_1,
..., \delta_{\lfloor \tfrac{g}{2}\rfloor})$ for $N^1(\M_g)$.  The following calculations still hold for $g=2$, but
$\check \lambda$ should be replaced by
$\tfrac{1}{10}\check\delta_{irr} + \tfrac{1}{5}\check \delta_1$ and
eliminated as a basis element.

Let $C_k \in \overline{NE}_1(X_{n,0})$ be as described in
\S\ref{S:effconeXn1} (define $C_1 \in \overline{NE}_1(X_{n,0})$ as the
image of $C_1\in \overline{NE}_1(X_{n,1})$).  Defining $\check
\delta_0 := 0 \in N^1(\M_g)$, the intersection numbers with
$i_*C_{2j+1}$ are 
\[
	\delta_{irr} = 2(2g+1-2j), \qquad 
	\delta_j = \frac{1}{2}(2j+1-2g), \qquad 
	\kappa + \delta_j = 2(g-j-1).
\]
Using the relation $12 \lambda \equiv \kappa + \delta$ gives 
\begin{align*}
	i_* C_{2j+1} \equiv 2(2g+1-2j)\check \delta_{irr} +
\frac{1}{2} (2j+1-2g)\check \delta_j + \frac{1}{2}(g-j)\check \lambda
\qquad (j \ge 0).
\end{align*}
Similarly, 
\[
	i_*C_{2j} \equiv 4(j-g)\check \delta_{irr} + (g+1-j) \check
\delta_j \qquad (j \ge 1). 
\]
By induction, 
\begin{equation}\label{E:images}  
\begin{aligned}
i_*\check b_{2j} &\equiv 2\check \delta_{irr} +
\frac{j(g+1-j)}{4g+2}\check \lambda \qquad (1 \le j \le g),\\
i_* \check b_{2j+1} &\equiv \frac{1}{2}\check \delta_j +
\frac{j(g-j)}{4g+2}\check \lambda \qquad (1 \le j \le g-1), 
\end{aligned}
\end{equation} 
where $\delta_j := \delta_{g-j}$ for $j > \lfloor
\tfrac{g}{2}\rfloor$. (These can also be deduced using the formulae
given in \cite[p. 303]{Mgbook}, but note the ramification of $X_n\to
X_{n,0}$ induced by $\s_{2g+2}$ is ignored in that calculation; the
most economical way to account for the ramification is to replace
their $\Xi_1$ with $\xi_1 := \tfrac{1}{2}\Xi_1$.)

\begin{corollary}[$g \ge 2$]
An effective divisor in $\M_g$ pulls back to $Nem(X_{2g+2,0})$ iff it
is contained in the cone generated by the following inequalities $(1
\le i \le g-1)$:
\begin{align*}
(2i-1)\check \delta_i \: &\le \: 4(2i+1)\check \delta_{irr} \: + \:
i\check \lambda,\\
4i \check \delta_{irr} \: &\le \: (i+1)\check \delta_i.
\end{align*}
\end{corollary}

\begin{proof}
This follows from (\ref{E:images}) and Theorem~\ref{T:NemZn}.  
\end{proof}

\subsection{The morphisms $X_{2n+3,1} \to \M_g, \: \M_{g,1}$ ($1 \le n
\le g-1,\,g$)}\label{SS:Mg1} To each element $C \in X_{2n+3,1}$ there
is a unique degree two admissible cover $\til C$ obtained as in
\S\ref{SS:Mg} via ignoring the distinguished point of $C$.  This
distinguished point then determines (up to hyperelliptic involution of
some component of $\til C$) a marked point $P \in \til C$.  Consider
the morphism $p:  X_{2n+3,1} \to \M_g$ ($1 \le n \le g-1$) given by
attaching to $\til C$ a fixed curve $D$ of complementary genus, by
identifying a fixed point $Q \in D$ with $P$.  Via a calculation
analogous to the one of the previous subsection,
\begin{align*}
	p_* C_1 &\equiv -(n-1)\check \delta_{g-n},\\
p_* C_{2j+1} &\equiv -4(n-j)\check \delta_{irr} +
(n+1-j)\check \delta_{g-n+j}, \quad \text{and}\\
p_* C _{2j} &\equiv 2(2n+3-2j) \check \delta_{irr}
-\frac{1}{2}(2n+1-2j) \check \delta_{g-n+j-1} +
\frac{1}{2}(n+1-j)\check \lambda.  
\end{align*}
\begin{remark}
As the above equations are written, they also apply when the target
space is $\M_{g,1}$; the morphism $p$ in that case is given by taking
some fixed point $(R \ne Q) \in D$ as marked point.
\end{remark}

Induction shows:
\begin{equation}\label{E:maptoMg1}
\begin{aligned}
p_*\check b_{2j+2} &\equiv \frac{1}{2}\check \delta_{g-n+j} +
\frac{j(n-j)}{2(2n+1)} \check \lambda 
- \frac{(n-j)(2n+1-2j)}{(2n+1)(n+1)} \check \delta_{g-n}, \\
p_* \check b_{2j+1} &\equiv 2\check \delta_{irr} +
\frac{j(n+1-j)}{2(2n+1)}\check \lambda -
\frac{(2n+1-2j)(n+1-j)}{(2n+1)(n+1)} \check \delta_{g-n}. \qquad 
\end{aligned}
\end{equation}
These equations also apply for $n = g$, if $\check \delta_0 := -
\check \omega \in N^1(\M_{g,1})$.

\begin{corollary}[$g \ge 3$, and $2\le n \le g-1$]\label{C:Mg1ineqs} 
The subcone of effective divisors of $\M_g$ pulling back to
$Nem(X_{2n+3,1})$ is the subset determined by $(1 \le k\le
n-1$ and $0\le m \le k-1$, unless specified otherwise$)$:
\begin{align*}
\bullet\: (k+1)\check \delta_{g-k} &\ge 4k\check \delta_{irr},\\
\bullet\: 4(2k+1)\check \delta_{irr} + k\check \lambda &\ge (2k-1)\check
\delta_{g-k},\\
\bullet\: (2m+1)(k-m)\check \delta_{g-k} + km(k-m)\check \lambda & 
+ k(2k+1)\check \delta_{g-n+m}\\ \ge k(2m+1)&\check \delta_{g-n+k}
+ 4k(k-m)\check \delta_{g-n},\\
\bullet\: m(k+1)(k-m)\check \lambda + (k+1)(2k+1)\check \delta_{g-n+m} 
&\ge 4m(2m+1) \check \delta_{irr}\\ 
 + (2k+1)&[2(k-m)+1]\check \delta_{g-n} \qquad  (0 \le m \le k),\\
\bullet\: (m+1)[2(k-m)-1]\check \delta_{g-k} + k(m+1)&(k-m)\check
\lambda + 4k(2k+1)\check \delta_{irr} \\ 
\ge 2k[2(k-m)-1]&\check \delta_{g-n} + 2k(m+1)\check \delta_{g-n+k}.  
\end{align*}
\end{corollary}

\begin{remark}
The corollary also applies (with the same notation) when the target
space is $\M_{g,1}$.  In that case, $g \ge 2$, $2\le n \le g$, and
$\check \delta_0 := - \check \omega$.  For $g=2$, $\check \lambda$
should be replaced by $\tfrac{1}{10}\check \delta_{irr} +
\tfrac{1}{5}\check \delta_1$.
\end{remark}

\begin{proof}
Pushing forward the inequalities defining $Nem(X_{2n+3,1})$
(Corollary~\ref{C:NemXn1}) shows the subcone in question is determined
by ($1 \le k\le n-1$ and $0\le m \le k-1$, unless specified
otherwise):
\begin{align*}
\mathbf{(i)}\: (k+1)\check \delta_{g-k} &\ge 4k\check \delta_{irr},\\
\mathbf{(ii)}\: 4(2k+1)\check \delta_{irr} + k\check \lambda \ge (2k-1)&\check
\delta_{g-k} \qquad (1\le k \le n),\\
\mathbf{(I)}\: (2m+1)(k-m)\check \delta_{g-k} + km(k-m)\check \lambda
&+ k(2k+1)\check \delta_{g-n+m}\\
 \ge k(2m+1)\check \delta_{g-n+k} &+ 4k(k-m)\check \delta_{g-n},\\
\mathbf{(II)}\: 4[2(k+m)+3]\check \delta_{irr} + (k+1)(m+1)\check \lambda &\ge 
2(2k+1)\check \delta_{g-n},\\
\mathbf{(III)}\: m(k+1)(k-m)\check \lambda + (k+1)(2k+1)\check
\delta_{g-n+m} &\ge 4m(2m+1) \check \delta_{irr} \\ +
(2k+1)[2(k-m)+1]\check \delta_{g-n} \qquad &(0 \le m \le k, \: 0 \le k
\le n-1),\\
\mathbf{(IV)}\: (m+1)[2(k-m)-1]\check \delta_{g-k} + k(m+1)(k-m)\check
\lambda &+ 4k(2k+1)\check \delta_{irr} \\
\ge 2k[2(k-m)-1]\check \delta_{g-n} &+ 2k(m+1)\check \delta_{g-n+k} . 
\end{align*}
Fix a positive integer $n\ge
2$.  For $n=2$ set $c_1 := 1$ and $c_2 := 2$.  For $n \ge 3$ set
\begin{align*}
c_1 &:= 2n(n-1)[2(k+m)+3] - 2(4n-3)(k+1)(m+1),\\
c_2 &:= 10(k+1)(m+1)-4[2(k+m)+3].  
\end{align*}
Then the $c_i$ are non-negative for $1 \le k \le n-1$, $0 \le m \le
k-1$, and (denoting e.g. inequality $(\mathbf{i})$ evaluated at $k=1$
by $(\mathbf{i})(1)$)
\[
c_1[(\mathbf{i})(1)+2(\mathbf{ii})(1)] +
c_2[(2n-3)(\mathbf{i})(n-1)+n(\mathbf{ii})(n-1)]
\]
gives
\[
	2(5n^2-13n+6)\cdot \biggl[ 4[2(k+m)+3]\check \delta_{irr}
+(k+1)(m+1)\check \lambda \biggr] \ge 0. 
\]
$5n^2-13n+6 > 0$ for $n \ge 3$, so ($3 \le n$, $1 \le k \le n-1$, $0
\le m \le k-1$):
\[
	4[2(k+m)+3]\check \delta_{irr} + (k+1)(m+1)\check \lambda
\ge 0.   
\]
Note that ``$\check \delta_{g-n} \le 0$'' comes from from
$\mathbf{III}(1,0)$, so $\mathbf{(II)}$, $\mathbf{III}(0,0)$, and
$\mathbf{(ii)}(n)$ are subsumed.
\end{proof} 

\section{The cone of effective divisors of $\M_{2,1}$}\label{S:M21} 
As an example we analyze $\n(\M_{2,1})$ in terms of $\n(X_{7,1})$.
Let
\[
	p: X_{7,1} \to \M_{2,1}
\]
be the morphism described in the first Remark of \S\ref{SS:Mg1}, for $n=g=2$; $p$ is easily seen to
be birational.  We continue to restrict to characteristic $\ne 2$ in this section.  

Let $W \subseteq \M_{2,1}$ denote the Weierstrass divisor, i.e. the
closure of the locus $\{(C,P) : P\in C$ is a Weierstrass point\}.  The
images of the boundary divisors of $X_{7,1}$ are
\[
\begin{aligned}
B_2 &\mapsto W,\\
B_3 &\mapsto \{\text{banana curves}\},
\end{aligned}
\qquad \qquad
\begin{aligned}
B_4 &\mapsto \Delta_{1;\{1\}},\\
B_5 &\mapsto \Delta_{irr}. \:
\end{aligned}
\]
The general element of the locus of ``banana curves" is a
genus one curve attached at two points to a smooth rational curve with
marked point.

From Proposition~\ref{P:effconeXn1}, $\overline{NE}^1(X_{7,1})
= \sum_{i=2}^5 \r B_i$.  
\begin{corollary}\label{C:M2,1}
$\n(\M_{2,1})$ is simplicial, generated by
the classes of $\Delta_{irr}$, $\Delta_1$ and $W$.  
\end{corollary}

\begin{proof}
The Picard number of $\M_{2,1}$ is three, and $B_3$ is contracted by
$p$.
\end{proof}

In terms of the ordered basis $(\Delta_{irr}, \Delta_{1}, W)$ it is
easy to check (using the $n=7$ example of \S\ref{S:NemXn1}) that the
pushforward of $Nem(X_{7,1})$ to $\M_{2,1}$ is generated by $A := (1,
1, 0)$, $B := (1, 6, 0)$, $D:= (1, 6, 20)$, and $E := (3, 3, 10)$. In
fact,

\begin{proposition}\label{P:MovM21}
$Mov(\M_{2,1}) = p_* Nem(X_{7,1}) = \r A + \r B + \r D + \r E$.
\end{proposition}

\begin{proof} 
By Proposition~\ref{P:functorial}, it suffices to show $A$, $B$, $D$,
and $E$ are generated by moving linear systems.  From
Lemma~\ref{L:f_*Mmoves} and Remark (b) following the $n=7$ example of
\S\ref{S:NemXn1}, the divisor $D$ moves, since it has a moving preimage
in $X_{7,1}$.  The classes $A$ and $B$ correspond to the classes
$\lambda$ and $12\lambda - \delta_{irr}$ of $\M_2$, respectively, so
both are basepoint free (the characteristic $p$
case follows from the results of \cite{nefbig}).  $E$ moves since it
can be written as an effective combination of either $W$ and $A$ or
$\Delta_{irr}$ and $D$.
\end{proof}

\begin{remarks}
\begin{enumerate}

\item $Nem(\M_{2,1})$ is also the cone of Proposition~\ref{P:MovM21},
as can be directly calculated \cite{Ru01}.


\item The pushforward of $Nef^1(X_{7,1})$ (presented in the $n=7$
example of \S\ref{S:NemXn1}) is the convex hull of the rays generated
by $A$, $B$, and $D$.

\item The equations (\ref{E:maptoMg1}) give
\[
	30 p^* \omega \equiv (5,12,6,2).  
\]
Recalling that $b_5 := \frac{1}{2}B_5$ and using the fact that 
$p$ has connected fibres, 
\[
	C := 30 \omega = 30 p_*p^*\omega \equiv \left( \Delta_{irr}
+ 6 \Delta_1 + 5 W\right).
\]
$Nef^1(\M_{2,1}) = \r A + \r B + \r C$ is easily calculable using test
curves and the facts that $A$, $B$, and $C$ are nef.
($Nef^1(\M_{2,1})$ is furthermore the pull-back of $Nef^1(\M_3)$
under the morphism induced by $\Delta_1 \subseteq \M_3$ \cite{Ru01}.)

\item The class $\omega$ of $\M_{2,1}$ is semi-ample in any
characteristic (the $\CHAR p$ case follows from the results of
\cite{nefbig}; the $\CHAR 0$ case is shown in \cite{Ru01}),
so every line bundle in $Nef^1(\M_{2,1})$ is eventually free.  

\item 
The inequalities of Corollary~\ref{C:Mg1ineqs} cut out the nef cone of
$\M_{2,1}$; i.e. any effective divisor on $\M_{2,1}$ pulling back to a
member of $Nem(X_{7,1})$ is nef.

\end{enumerate}
\end{remarks}

\begin{figure}[!ht]
\begin{center}
	\scalebox{.35}{\includegraphics{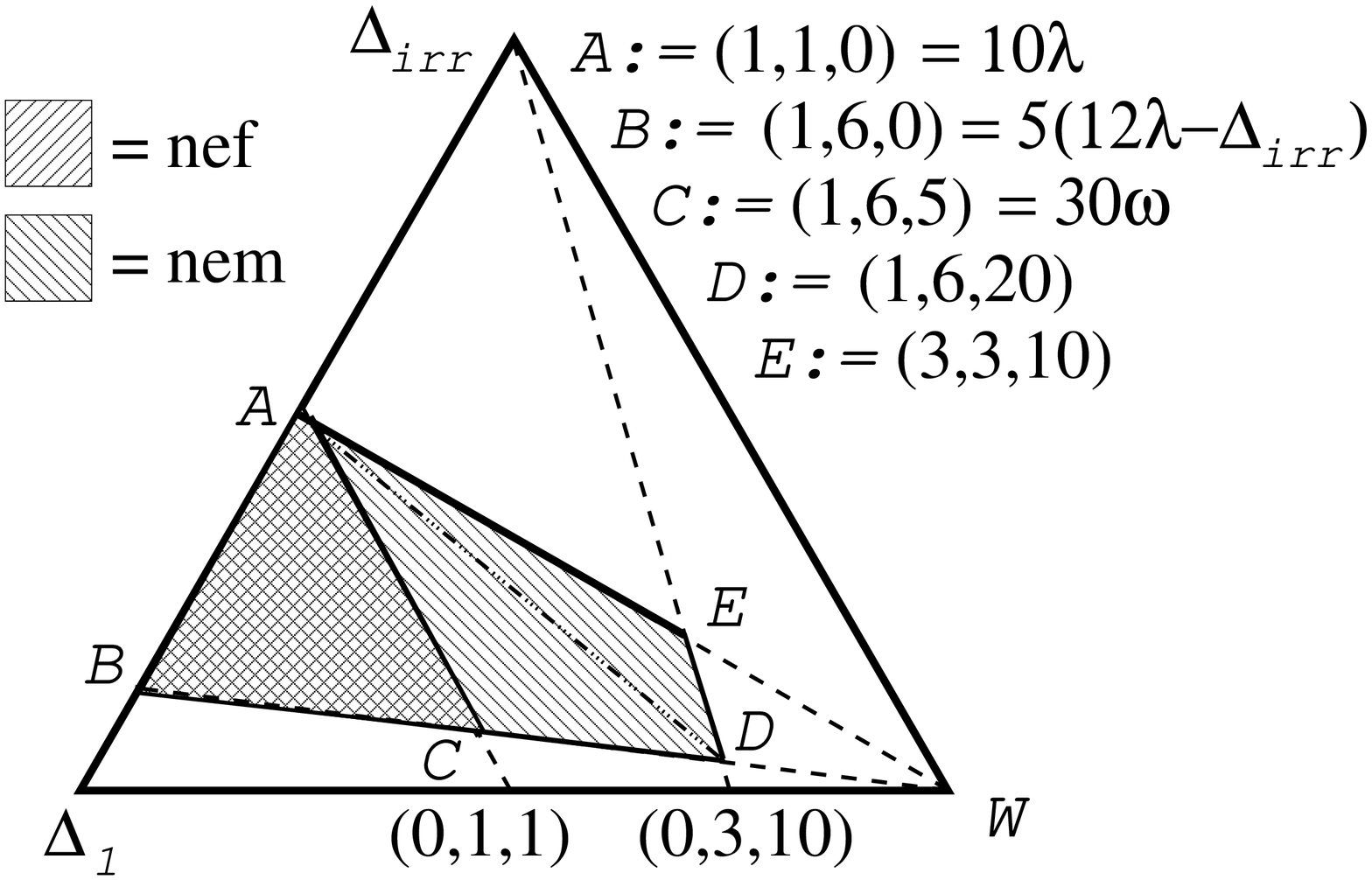}} 
	\caption{A cross section of 
	$\n(\M_{2,1})$.}\label{Fi:coneM2,1}\label{coneM2,1}
\end{center}
\end{figure}

We partially decompose $Mov(\M_{2,1})$, in characteristic zero.  
\begin{proposition}[$\Char 0$]
$X_{7,1}$ resolves a rational map $\M_{2,1} \dashrightarrow
Z$, where $Z$ is a $\Q$-factorial projective variety isomorphic to
$\M_{2,1}$ in codimension one.  The nef cone of $Z$ is identified with
the convex hull of the divisors $A$, $C$, and $D$ (Figure~\ref{Fi:coneM2,1}).
\end{proposition}
\begin{proof}

The morphism $p: X_{7,1} \to \M_{2,1}$ contracts the curve class 
\[
C_1 \equiv \check b_5 + 2\check b_2 - \check b_3,
\]
which is obtained from a general member of $B_3$ by varying the node
along the component containing the distinguished marked point.  Since
the relative Picard number of $p$ is one, $C_1$ is the only curve
class contracted by $p$; since $B_3$ intersects $C_1$ negatively, any
curve numerically equivalent to $C_1$ is contained in $B_3$, so $p$ is
an isomorphism outside of $B_3$.

Let $C_2$ be the curve class of $X_{7,1}$ obtained by varying the node
of a general member of $B_3$ along the other component, so
\[
	C_2 \equiv 4\check b_4 - 2\check b_3.  
\]
Intersection with $C_2$ determines a supporting hyperplane of
$Nef^1(X_{7,1})$, and $C_2^{\perp}\cap Nef^1(X_{7,1})$ is the convex
hull of $(0,2,1,2), (5,12,6,2)$, and $(10,6,3,1)$, so is a
full-dimensional face; thus $C_2$ is an extremal ray of the Mori cone
$\overline{NE}_1(X_{7,1})$.  We will use the Cone Theorem
\cite[Theorem 3.7, p. 76]{blue} to get a morphism $q$ of $X_{7,1}$
contracting only the class $C_2$.

For $n=7$ Riemann-Hurwitz and \cite[Lemma 3.5]{KeMc96} or
\cite[Prop. 1]{Pand} give 
\[
	K_{X_{7,1}} \equiv -\tfrac{1}{3}b_2 - \tfrac{4}{3} b_5, 
\]
which dots $C_2$ trivially.  $X_{n,1}$ is a finite quotient of the
smooth variety $X_n$, and the boundary of $X_n$ is snc, so $(X_{7,1},
\frac{1}{2}B_3)$ is klt.  By the Cone Theorem there exists a morphism
$q: X_{7,1} \to Z$ with connected fibres onto a normal, projective
variety $Z$, contracting only curves numerically equivalent to $C_2$
(so $q$ contracts $B_3$ and is an isomorphism outside of it).  Since
the relative Picard number of $q$ is one and $B_3$ is contracted, $Z$
is $\Q$-factorial.

Finally, $q^*Nef^1(Z)$ is the face of $Nef^1(X_{7,1})$ determined by
$C_2^{\perp}$, so $p_* q^* Nef^1(Z)$ is the convex hull of the rays
generated by $A$, $C$, and $D$.
\end{proof} 

%
%
\bibliographystyle{amsalpha}
\bibliography{bib}
\end{document}